\documentclass[leqno,12pt]{amsart} %leqno is the option to put formula numbers on the left side
\usepackage{amsmath, amssymb}
\usepackage{amsthm} %theorem environment option
\usepackage{amscd,graphicx,
latexsym,amsthm,verbatim}
\usepackage[mathscr]{eucal}
\usepackage[all,cmtip]{xy}
\usepackage{pb-diagram}

\theoremstyle{plain} %text of this environment is typesetted in italics
\newtheorem{theorem}{\indent\sc Theorem}[section]
\newtheorem{lemma}[theorem]{\indent\sc Lemma}
\newtheorem{corollary}[theorem]{\indent\sc Corollary}

\theoremstyle{definition} %text of this environment is typesetted in roman letters
\newtheorem{definition}[theorem]{\indent\sc Definition}
\newtheorem{remark}[theorem]{\indent\sc Remark}

\title{{Projectivized rank two toric vector bundles are mori dream spaces}} %title of the paper
\author[Jos\'e Luis Gonz\'alez]{
{Jos\'e Luis Gonz\'alez$^{*}$} %names of authors
}
\date{} %leave empty
\address{Jos\'e Luis Gonz\'alez, Department of Mathematics, University of
Michigan.}
\email{jgonza@umich.edu}
\newcommand{\pe}{{{\bf{P}}}(\mathcal{E})}
\newcommand{\Ox}{\mathcal{O}_{X}}
\newcommand{\Ope}{\mathcal{O}_{{\bf{P}}(\mathcal{E})}}
\newcommand{\LL}{\mathcal{L}}
\newcommand{\EE}{\mathcal{E}}
\newcommand{\om}{\underline{{\bf{m}}}}
\newcommand{\omp}{\underline{{\bf{m}}}'}
\newcommand*{\longhookrightarrow}{\ensuremath{\lhook\joinrel\relbar\joinrel\rightarrow}} 
\newcommand{\sbl}{\vskip 3pt}

\begin{document}

%%%%%%%%%%%%%%% footnote %%%%%%%%%%%%%%%%
\footnote{ %2000 MSC numbers
}
\footnote{ %key words and phrases
}
\footnote{ %acknowledgment of support etc. if any
$^{*}$The author was partially supported by the NSF under grant DMS-0502170 and by the David and Lucile Packard Foundation under M. Musta\c{t}\u{a}'s fellowship.
}
%%%%%%%%%%%%%%%%%%%%%%%%%%%%%%%%%%%%%%%%%

\begin{abstract}
We prove that the Cox ring of the projectivization ${{\bf{P}}}(\mathcal{E})$ of a rank two toric vector bundle 
$\mathcal{E}$, over a toric variety $X$, is a finitely generated $k$-algebra. As a consequence, ${{\bf{P}}}(\mathcal{E})$ is a Mori dream space if the toric variety $X$ is projective and simplicial. %In the proof, we describe a set of generators of the Cox ring of ${{\bf{P}}}(\mathcal{E})$ in terms of data that arises from the Klyachko filtrations of $\mathcal{E}$.
\end{abstract}

\maketitle

%\tableofcontents

%%%%%%%%%%%%%%%%%%%%%%%%%%%%%%%%%%%%%%%%%%%%%%%%%%
%%                 INTRODUCTION                 %%
%%%%%%%%%%%%%%%%%%%%%%%%%%%%%%%%%%%%%%%%%%%%%%%%%%

\section{Introduction}
 
%\vspace{5 mm}
Mori dream spaces were introduced by Hu and Keel in \cite{HK} as a class of varieties with interesting features from the point of view of Mori theory. For instance, their pseudoeffective and nef cones are both polyhedral, and the Mori program can be carried out for any pseudoeffective divisor on these varieties. Additionally, their pseudoeffective cones can be decomposed into finitely many closed convex chambers that are in correspondence with the birational contractions of $X$ having ${\bf{Q}}$-factorial image (see Proposition 1.11 \cite{HK}).  
A Mori dream space can be defined as a normal projective ${\bf{Q}}$-factorial variety $X$, with $\operatorname{Pic}(X)_{{\bf{Q}}} = N^{1}(X)_{{\bf{Q}}}$ and a finitely generated Cox ring (see Definition \ref{definition.cox.ring}). One basic example is that of toric varieties, where the Cox ring is a polynomial ring in finitely many variables (see \cite{Cox}). Thus, a projective simplicial toric variety is a Mori dream space. One more example is given by the (log) Fano varieties, which have recently been proven to be Mori dream spaces (see \cite{BCHM}). Other references for recent work on Cox rings of particular Mori dream spaces are \cite{BP}, \cite{Cas} and \cite{STV}.   

In \cite{PTVB}, Hering, Musta\c{t}\u{a} and Payne ask whether projectivizations of toric vector bundles are Mori dream spaces, and give some evidence for a positive answer. For example, the Mori cones of these varieties are rational polyhedral, and they are known to be Mori dream spaces for toric vector bundles that split as a sum of line bundles. Toric vector bundles were classified by Klyachko (see \cite{Klyachko}) via an equivalence of categories under which each toric vector bundle corresponds to a collection of filtrations of a suitable vector space (see Theorem \ref{Theorem.Klyachko}).

In this paper, we prove that the Cox ring of the projectivization $\pe$ of a rank two toric vector bundle $\EE$ over an arbitrary toric variety $X$ is finitely generated as a $k$-algebra. As a consequence, $\pe$ is a Mori dream space, if $X$ is assumed to be projective and simplicial. In the proof, we consider a particular set of generators of the group $N^{1}(\pe)$, and describe a finite set of generators for the associated Cox ring. This is done by considering a finer grading of the Cox ring induced by the torus action, and describing the graded pieces and the multiplication map in terms of data arising from Klyachko's filtrations of $\EE$.

Alternative proofs of our main result can be found in \cite{Knop} and \cite{Suess}. In fact, in these articles the authors establish the finite generation of Cox rings in the more general setting of $T$-varieties of complexity one. We recall that a normal variety with an algebraic action of a torus $T$ is called a $T$-variety, and that its complexity is defined to be the lowest codimension of a $T$-orbit on the variety. Hence projectivizations of rank two toric vector bundles are examples of $T$-varieties of complexity one. On the other hand, we hope that our methods will provide new ideas for treating the analogous problem in the case of projectivizations of higher rank toric vector bundles.

\subsection*{Acknowledgments}
I am grateful to Milena Hering, Mircea Musta\c{t}\u{a} and Sam Payne from whom I learned the questions that motivated this work. I thank M. Hering for bringing the paper \cite{Knop} to my attention, and I especially thank Mircea Musta\c{t}\u{a} for many inspiring conversations.

%I am grateful to Milena Hering, Mircea Musta\c{t}\u{a} and Sam Payne from whom I learned the questions that motivated this work. I thank M. Hering for communicating the results in the paper \cite{Knop} to me, and I especially thank Mircea Musta\c{t}\u{a} for many inspiring conversations.

%%%%%%%%%%%%%%%%%%%%%%%%%%%%%%%%%%%%%%%%%%%%%%%%%%%%%%%%%%%%%%%%%%%%%%%%%%%%%%%%%%%%%%%%%

%\vspace{6 mm}

%%%%%%%%%%%%%%%%%%%%%%%%%%%%%%%%%%%%%%%%%%%%%%%%%%%%%%%%%%%%
%%           Toric Vector Bundles and Cox rings           %%  
%%%%%%%%%%%%%%%%%%%%%%%%%%%%%%%%%%%%%%%%%%%%%%%%%%%%%%%%%%%%

\section{Toric Vector Bundles and Cox rings} \label{TVBCR}

All our varieties are defined over a fixed algebraically closed field $k$. By a \emph{line bundle} on a variety $Z$, we mean a locally free sheaf of rank one on $Z$. The groups of line bundles on $Z$ modulo linear equivalence and modulo numerical equivalence are denoted by $\operatorname{Pic}(Z)$ and $N^{1}(Z)$, and their ${\bf{Q}}$-extensions of scalars are denoted by $\operatorname{Pic}(Z)_{{\bf{Q}}}$ and $N^{1}(Z)_{{\bf{Q}}}$, respectively. We follow the convention that the \emph{geometric vector bundle} associated to the locally free sheaf $\mathcal{E}$ is the variety ${\bf{V}}(\mathcal{E}) =  \text{\textbf{Spec}} \, \bigoplus_{m \geq 0} Sym^{m}\mathcal{E}^{\vee}$, whose sheaf of sections is $\mathcal{E}$. Also, when we consider the fiber of $\mathcal{E}$ over a point $z \in Z$, we mean the fiber over $z$ of the projection $f \colon {\bf{V}}(\mathcal{E}) \rightarrow Z$. Lastly, by the \emph{projectivization} ${\bf{P}}(\mathcal{E})$ of $\mathcal{E}$, we mean the projective bundle $\text{\textbf{Proj}} \, \bigoplus_{m \geq 0} Sym^{m}\mathcal{E}$ over $Z$. This bundle is endowed with a projection $\pi \colon {\bf{P}}(\mathcal{E}) \rightarrow Z $ and an invertible sheaf $\mathcal{O}_{{\bf{P}}(\mathcal{E})}(1)$ (see II.7 in \cite{Hartshorne}).

Let $X$ be an $n$-dimensional toric variety corresponding to the nondegenerate fan $\Delta$ in the lattice $N \cong {\bf{Z}}^{n}$. We denote the algebraic torus acting on $X$ by $T$, and we designate its unit element by $t_{0}$. We denote the character lattice $\operatorname{Hom}(N,{\bf{Z}}) \cong {\bf{Z}}^{n}$ of $T$ by $M$. Thus, $T=\operatorname{Spec}k[M] = \operatorname{Spec}k[ \, \chi^{u} \, | \, u \in M \, ]$ and $X$ has an open covering given by the affine toric varieties $U_{\sigma} = \operatorname{Spec}k[ \, \sigma^{\vee} \cap M \, ]$ corresponding to the cones $\sigma \in \Delta$. We denote the rays in $\Delta$ by $\rho_{1},\rho_{2},\ldots,\rho_{d}$. For each ray $\rho_{j}$, we denote its primitive lattice generator by $v_{j}$ and its associated torus invariant prime divisor by $D_{j}$. For a detailed treatment of toric varieties we refer to \cite{Fulton}.

% Definition of toric vector bundle

A \emph{toric vector bundle} on the toric variety $X$ is a locally free sheaf $\mathcal{E}$ together with an action of the torus $T$ on the variety ${\bf{V}}(\mathcal{E})$, such that the projection $f \colon {\bf{V}}(\mathcal{E}) \rightarrow X$  is equivariant and $T$ acts linearly on the fibers of $f$. Given any $T$-invariant open subset $U$ of $X$, there is an induced action of $T$ on $H^{0}(U,\mathcal{E})$, defined by the equation
\[
(t \cdot s)(x)=_{\operatorname{def}} t (s( t^{-1} x )),
\]
for any $s \in H^{0}(U,\mathcal{E})$, $t \in T$ and $x \in X$. This action induces a direct sum decomposition 
\begin{equation} \label{isotypical}
H^{0}(U,\mathcal{E}) = \bigoplus_{u \in M} H^{0}(U,\mathcal{E})_{u},
\end{equation}
where $H^{0}(U,\mathcal{E})_{u}=_{def}\{s \in H^{0}(U,\mathcal{E}) \ | \ t \cdot s = \chi^{u} (t) s \textnormal{ for each $t \in T$}  \}$.
We refer to the summands in the decomposition from (\ref{isotypical}) as the isotypical components of $H^{0}(U,\mathcal{E})$ with respect to the torus action. 
Let us denote the fiber of $\EE$ over $t_{0}$ by $E$, and for any $T$-invariant open subset $U$, let $ev_{t_{0}} \colon H^{0}(U,\mathcal{E}) \rightarrow E$ be the evaluation map at $t_{0}$. For each ray $\rho_{j} \in \Delta$ and each $u \in M$, the evaluation map $ev_{t_{0}}$ gives an inclusion $H^{0}(U_{\rho_{j}},\mathcal{E})_{u} \longhookrightarrow E$. If $u, u' \in M$ satisfy $\langle u ,v_{j} \rangle \geq \langle u' ,v_{j} \rangle$, then $ev_{t_{0}}(H^{0}(U_{\rho},\mathcal{E})_{u}) \subseteq ev_{t_{0}}(H^{0}(U_{\rho},\mathcal{E})_{u'})$. In particular, $ev_{t_{0}}(H^{0}(U_{\rho},\mathcal{E})_{u})$ depends only on $\langle u ,v_{j} \rangle$, or equivalently, only on the class of $u$ in $M/\rho_{j}^{\perp}\cap M \cong {\bf{Z}}$. Hence we may simply denote this space by $\mathcal{E}^{\rho_{j}}(\langle u ,v_{j} \rangle)$.
The ordered collection $\EE^{\rho_{j}} = \{ \EE^{\rho_{j}}(i) \ | \  i \in {\bf{Z}} \}$ gives a decreasing filtration of $E$, and the set of these filtrations $\{ \EE^{\rho_{j}} \ | \ 1 \leq j \leq d \}$ are known as the \emph{Klyachko filtrations} associated to $\EE$. Klyachko proved in \cite{Klyachko} that these filtrations of $E$ satisfy a certain compatibility condition, and that they completely describe $\EE$. More precisely,

\begin{theorem}[Klyachko] \label{Theorem.Klyachko} The category of toric vector bundles on the toric variety $X$ is e\-qui\-va\-lent to the category of finite dimensional $k$-vector spaces $E$ with collections of decreasing filtrations $\{\mathcal{E}^{\rho}(i) \ | \ i \in {\bf{Z}} \}$, indexed by the rays of $\Delta$, satisfying the following compatibility condition: For each cone $\sigma \in \Delta$, there is a decomposition $E= \bigoplus_{\bar{u} \in M/\sigma^{\perp} \cap M}E_{\bar{u}}$ such that
\[
\mathcal{E}^{\rho}(i)=\sum_{\langle \bar{u},v_{\rho} \rangle \geq i}{E_{\bar{u}}},
\]
for every ray $\rho \subseteq \sigma$ and every $i \in {\bf{Z}}$.
\end{theorem}

Some references for recent work on toric vector bundles are \cite{Jose1}, \cite{PTVB}, \cite{Payne2} and \cite{Payne}. Next, we give the definition of the algebras that we will call Cox rings of varieties, and we briefly compare this definition with the one used in \cite{HK}. 

\begin{definition} \label{definition.cox.ring}
Let $X$ be a variety such that $\operatorname{Pic}(X)_{\bf{Q}}$ is finite dimensional. Any $k$-algebra of the form
\[
\operatorname{Cox}\bigl(X,(\LL_{1},\LL_{2},\ldots,\LL_{s})\bigr)=_{def} \bigoplus_{(m_{1},m_{2},\ldots,m_{s}) \in {\bf{Z}}^{s}} H^{0}(X,\LL_{1}^{\otimes m_{1}} \otimes \LL_{2}^{\otimes m_{2}} \otimes \cdots \otimes \LL_{s}^{\otimes m_{s}})
\]
where $\LL_{1},\LL_{2},\ldots,\LL_{s}$ are line bundles on $X$ whose classes span $\operatorname{Pic}(X)_{{\bf{Q}}}$, will be called a \emph{Cox ring} of $X$.
\end{definition}

In the setting of Definition \ref{definition.cox.ring}, it is easy to see that the finite generation of one Cox ring of $X$ is equivalent to the finite generation of every Cox ring of $X$. In \cite{HK}, Hu and Keel give a definition of Cox rings for any projective variety $X$ that satisfies $\operatorname{Pic}(X)_{{\bf{Q}}} = N^{1}(X)_{{\bf{Q}}}$. Their definition is similar to Definition \ref{definition.cox.ring}, but they require the line bundles $\LL_{1},\LL_{2},\ldots,\LL_{s}$ to satisfy some extra conditions. Namely, the classes of $\LL_{1},\LL_{2},\ldots,\LL_{s}$ are required to form a basis for $\operatorname{Pic}(X)_{{\bf{Q}}}$ and their affine hull must contain the pseudoeffective cone of $X$. For projective varieties that satisfy $\operatorname{Pic}(X)_{{\bf{Q}}} = N^{1}(X)_{{\bf{Q}}}$, it is easy to see that a particular (equivalently every) Cox ring of $X$ in the sense of Hu and Keel is finitely generated if and only if a particular (equivalently every) Cox Ring of $X$ in the sense of Definition \ref{definition.cox.ring} is finitely generated. Since our interest lies in the finite generation of these $k$-algebras, and this does not depend on the definition of Cox ring that we use, we will use Definition \ref{definition.cox.ring} since it applies in the case of the projectivization of a toric vector bundle over an arbitrary toric variety.

%%%%%%%%%%%%%%%%%%%%%%%%%%%%%%%%%%%%%%%%%%%%%%%%%%%%%%%%%%%%%%%%%%%%%%%%%%%%%%%%%%%%%%%%%

% Reduction to the smooth case.

\begin{remark} \label{smooth}
Given a toric variety $X$ there exists a toric resolution of singularities, i.e. a smooth toric variety $X'$ and a proper birational toric morphism $f \colon X' \rightarrow X$. Given a toric vector bundle $\mathcal{E}$ on $X$, the induced map $f' \colon {\bf{P}}(f^{*}{\mathcal{E}}) \rightarrow \pe$ is also proper and birational. In this case $f^{*}{\mathcal{E}}$ is a toric vector bundle on $X'$ and the finite generation of a Cox ring of ${\bf{P}}(f^{*}{\mathcal{E}})$ implies the finite generation of any Cox ring of $\pe$. To see this we consider line bundles $\LL_{1},\LL_{2},\ldots,\LL_{s}$ on $\pe$ whose classes span $\operatorname{Pic}(\pe)_{{\bf{Q}}}$ and line bundles $\LL'_{1},\LL'_{2},\ldots,\LL'_{s'}$ on ${\bf{P}}(f^{*}{\mathcal{E}})$ whose classes span $\operatorname{Pic}({\bf{P}}(f^{*}{\mathcal{E}}))_{{\bf{Q}}}$. The finite generation of the algebra $\operatorname{Cox}\bigl({\bf{P}}(f^{*}{\mathcal{E}}),(\LL'_{1},\LL'_{2},\ldots,\LL'_{s'},f'^{*}\LL_{1},f'^{*}\LL_{2},\ldots,f'^{*}\LL_{s})\bigr)$ implies the finite generation of the subalgebra \[
\bigoplus_{(m_{1},m_{2},\ldots,m_{s}) \in {\bf{Z}}^{s}} \hspace{-0.3cm} H^{0}({\bf{P}}(f^{*}{\mathcal{E}}),f'^{*}\LL_{1}^{\otimes m_{1}} \otimes f'^{*}\LL_{2}^{\otimes m_{2}} \otimes \cdots \otimes f'^{*}\LL_{s}^{\otimes m_{s}})
= \operatorname{Cox}\bigl(\pe,(\LL_{1},\LL_{2},\ldots,\LL_{s})\bigr).
\] 
Therefore, to prove the finite generation of the Cox rings of projectivizations of rank two toric vector bundles over arbitrary toric varieties, it is enough to consider the case when the base is a smooth toric variety.
\end{remark}

%%%%%%%%%%%%%%%%%%%%%%%%%%%%%%%%%%%%%%%%%%%%%%%%%%%%%%%%%%%%%%%%%%%%%%%%%%%%%%%%%%%%%%%%%

% Presentation of the cox ring as R
Now, let $X$ be a smooth toric variety and let $\pi \colon \pe \rightarrow X$ be the projectivization of the rank two toric vector bundle $\EE$ over $X$. The classes of the line bundles $\Ox(D_{1}), \linebreak[0] \Ox(D_{2}), \linebreak[0] \ldots, \linebreak[0] \Ox(D_{d})$ span $\operatorname{Pic}(X)_{{\bf{Q}}} = N^{1}(X)_{{\bf{Q}}}$. Since $\operatorname{Pic}(\pe) = \operatorname{Pic}(X) \oplus {\bf{Z}} \cdot [\Ope(1)]$, the classes of the line bundles $\Ope(1), \linebreak[0] \pi^{*} \Ox(D_{1}), \linebreak[0] \pi^{*} \Ox(D_{2}), \linebreak[0] \ldots, \linebreak[0] \pi^{*} \Ox(D_{d})$ span $\operatorname{Pic}(\pe)_{{\bf{Q}}} = N^{1}(\pe)_{{\bf{Q}}}$. Therefore the following algebra is a Cox ring of $\pe$ in the sense of Definition \ref{definition.cox.ring}:
\[ 
%Cox(\pe, \Ope(1),\pi^{*} \Ox(D_{1}), \ldots,\pi^{*} \Ox(D_{d})) = 
C \ = \bigoplus_{(m,m_{1},\ldots,m_{d}) \in {\bf{Z}}^{d+1}} H^{0} (\pe, \Ope(m) \otimes \pi^{*} \Ox(m_{1}D_{1}) \otimes \cdots \otimes \pi^{*} \Ox(m_{d}D_{d}) ).
\]
By the projection formula, $C$ is isomorphic to the algebra
\[
R \ = \bigoplus_{(m,m_{1},\ldots,m_{d}) \in {\bf{Z}}^{d+1} } H^{0} (X, Sym^{m} \EE \otimes  \Ox(m_{1}D_{1}) \otimes \cdots \otimes \Ox(m_{d}D_{d})), 
\]
where $Sym^{m} \EE = 0$ for $m < 0$.
For each $\om = (m_{1},\ldots,m_{d}) \in {\bf{Z}}^{d}$, let us denote the toric line bundle 
$\Ox(m_{1}D_{1}) \otimes \cdots \otimes \Ox(m_{d}D_{d})$ on $X$ by $\LL^{\om}$, and the fiber of $\LL^{\om}$ over $t_{0}$ by $L_{\om}$.
% When X is not smooth:
%For each $m \in {\bf{Z}}$ and each $\om \in {\bf{Z}}^{d}$, we have an induced action of $T$ on $H^{0}(X,Sym^{m}\EE \otimes \LL^{\om})$, and for each $u \in M$ we define
%\[ 
%H^{0}(X,Sym^{m}\EE \otimes \LL^{\om})_{u} = \{s \in H^{0}(X,Sym^{m}\EE \otimes \LL^{\om}) \ | \ t \cdot s = \chi^{u} (t) s \textnormal{ for each $t \in T$}  \}
%\]
%When $\LL^{\om}$ is a line bundle, $Sym^{m}\EE \otimes \LL^{\om}$ is a toric vector bundle and this notation coincides with the one that we introduced before.
Note that for each $m,m' \in {\bf{Z}}$, each $\om, \omp \in {\bf{Z}}^{d}$ and each $u,u' \in M$, the product of $H^{0}(X,Sym^{m}\EE \otimes \LL^{\om})_{u}$ and $H^{0}(X,Sym^{m'}\EE \otimes \LL^{\omp})_{u'}$ in the algebra $R$ is contained in $H^{0}(X,Sym^{m + m'}\EE \otimes \LL^{\om+ \omp })_{u+u'}$.
Therefore we get a finer grading for $R$ given by
\[
R  
\ = \ \bigoplus_{(u,m,\om) \in M \times {\bf{Z}} \times {\bf{Z}}^{d}}  R_{(u,m,\om)}
\ = \ \bigoplus_{(u,m,\om) \in M \times {\bf{Z}} \times {\bf{Z}}^{d}}  H^{0}(X, Sym^{m}\EE \otimes \LL^{\om})_{u}. 
\]
Note that each of the homogeneous components $R_{(u,m,\om)} = H^{0}(X, Sym^{m}\EE \otimes \LL^{\om})_{u}$ is a finite dimensional $k$-vector space.
For each $l \in {\bf{Z}}^{+}$, let $R^{(l)}$ be the Veronese subalgebra of $R$ given by
\begin{equation} \label{veronese.cox.ring}
R^{(l)}
\ = \ \bigoplus_{(u,m,\om) \in M \times {\bf{Z}} \times {\bf{Z}}^{d}}  R_{(lu,lm,l\om)} 
\ = \ \bigoplus_{(u,m,\om) \in M \times {\bf{Z}} \times {\bf{Z}}^{d}}  H^{0}(X, Sym^{lm}\EE \otimes \LL^{l\om})_{lu}. 
\end{equation}
Since $R$ is a domain and $H^{0}(X,\Ox)$ is a finitely generated $k$-algebra, it follows from general considerations that the finite generation of $R$ is equivalent to the finite generation of $R^{(l)}$ for any $l \in {\bf{Z}}^{+}$.

%%%%%%%%%%%%%%%%%%%%%%%%%%%%%%%%%%%%%%%%%%%%%%%%%%%%%%%%%%%%%%%%%%%%%%%%%%%%%%%%%%%%%%%%%

%\vspace{6 mm}

%%%%%%%%%%%%%%%%%%%%%%%%%%%%%%%%%%%%
%%       Preliminary lemmas       %%  
%%%%%%%%%%%%%%%%%%%%%%%%%%%%%%%%%%%%

\section{Preliminary lemmas} \label{lemmas}

% Computations with symetric products

%\begin{notation} 
Let $W_{1},W_{2},\ldots,W_{l}$ be subspaces of a vector space W. For every collection of nonnegative integers $m, c_{1},c_{2},\ldots,c_{l}$, we denote by $Sym_{W}^{m}(W_{1}^{c_{1}} , W_{2}^{c_{2}} , \cdots , W_{l}^{c_{l}})$ the subspace of $Sym^{m}W$ equal to the image of the composition of the natural maps
	\[
W_{1}^{\otimes c_{1}} \otimes W_{2}^{\otimes c_{2}} \otimes \cdots \otimes W_{l}^{\otimes c_{l}} \otimes W^{\otimes (m - \sum_{i=1}^{l}{c_{i}})} \longrightarrow W^{\otimes m} \longrightarrow Sym^{m}W,
\]
if  $m \geq \sum_{i=1}^{l}{c_{i}}$, or the subspace $0$ of $Sym^{m}W$, otherwise. 
%\end{notation}

\begin{lemma} \label{lemma:M}
Let $W_{1},W_{2},\ldots,W_{q}$ be distinct subspaces of a vector space $W$. Let $m, m', \linebreak[0] c_{1}, \linebreak[0] c_{2}, \linebreak[0] \ldots , \linebreak[0] c_{q} , \linebreak[0] c_{1}' , c_{2}' , \ldots , c_{q}'$ be nonnegative integers.
 \item[{\bf{(a)}}] If $\sum_{i=1}^{q} c_{i} \leq m$ and $\sum_{i=1}^{q} c_{i}' \leq m'$, then
\[
\mu \biggl( Sym^{m}_{W}\bigl(W_{1}^{c_{1}},\ldots,W_{q}^{c_{q}} \bigr)
\otimes
Sym^{m'}_{W}\bigl(W_{1}^{c_{1}'},\ldots,W_{q}^{c_{q}'} \bigr) \biggr) = Sym^{m+m'}_{W}\bigl(W_{1}^{c_{1} + c_{1}'},\ldots,W_{q}^{c_{q} + c_{q}'} \bigr),
\]
where $\mu \colon Sym^{m}W \otimes Sym^{m'}W \rightarrow Sym^{m+m'}W$ is the multiplication map. 
 \item[{\bf{(b)}}] If $W_{1},W_{2},\ldots,W_{q}$ are one-dimensional, and $W$ is two-dimensional, then
\[
Sym_{W}^{m}(W_{1}^{c_{1}}) \cap Sym_{W}^{m}(W_{2}^{c_{2}}) \cap \cdots \cap Sym_{W}^{m}(W_{q}^{c_{q}}) = Sym_{W}^{m}(W_{1}^{c_{1}} , W_{2}^{c_{2}} , \cdots , W_{q}^{c_{q}}).
\]  
Furthermore, this subspace of $Sym^{m}W$ is nonzero precisely when $m \geq \sum^{q}_{i=1}{c_{i}}$, and in that case its dimension is $m+1-\sum^{q}_{i=1}{c_{i}}$.
\end{lemma}
\begin{proof}
 \item[{\bf{(a)}}] The conclusion follows at once from the commutativity of the diagram
\[
% \ar@{shape}[dir]  Shapes: |-> , {-}, {^{(}->}, etc. Directions: u,d,l,r.
\xymatrixcolsep{1.5pc}
%\xymatrixrowsep{2pc}
\xymatrix{
{\scriptstyle{ \bigl( ( \bigotimes\limits_{i} W_{i}^{\otimes c_{i}} ) \otimes W^{\otimes(m - \sum c_{i}) } \bigr)
\otimes
\bigl( ( \bigotimes\limits_{i} W_{i}^{\otimes c_{i}'} ) \otimes W^{\otimes(m' - \sum c_{i}') } \bigr) }}
\ar[d]^-{{\textstyle\wr}} \ar[r] &  {\scriptstyle{ W^{\otimes m} \otimes W^{\otimes m'} }} \ar[r] & {\scriptstyle{ Sym^{m}W \otimes Sym^{m'}W }} \ar[d]^-{\mu} \\
{\scriptstyle{  { ( \bigotimes\limits_{i} W_{i}^{\otimes ( c_{i} + c_{i}' ) } ) \otimes W^{\otimes( (m + m') - \sum (c_{i} + c_{i}' )) }  } }} 
\ar[r] & {\scriptstyle{ W^{\otimes (m+m')} }} \ar[r] & {\scriptstyle{ Sym^{m+m'}W.} }}
\]
 \item[{\bf{(b)}}] We fix an isomorphism of $k$-algebras between $\oplus_{h\geq 0} Sym^{h}W$ and the polynomial ring in two variables $k[x,y]$. The subspaces $W_{1},\ldots,W_{q}$ of $W$ correspond to the linear spans of some distinct linear forms $l_{1},\ldots,l_{q}$. The subspaces $ \cap_{i=1}^{q}Sym_{W}^{m}(W_{i}^{c_{i}})$ and $Sym_{W}^{m}(W_{1}^{c_{1}} , \linebreak[0] \cdots , \linebreak[0] W_{q}^{c_{q}})$ of $Sym^{m}W$ both correspond to the homogeneous polynomials of degree $m$ divisible by $l_{1}^{c_{1}}\cdots l_{q}^{c_{q}}$. From this observation the conclusion follows.
\end{proof}

%%%%%%%%%%%%%%%%%%%%%%%%%%%%%%%%%%%%%%%%%%%%%%%%%%%%%%%%%%%%%%%%%%%%%%%%%%%%%%%%%%%%%%%%%

% Definition of the sets A_{i}, the set \mathcal{J} and the linear functions $\lambda_{j}$.

Let $\EE$ be a toric vector bundle of rank two on $X$. Let $E$ be the fiber of $\EE$ over the unit of the torus, and let $\EE^{\rho_{1}}, \EE^{\rho_{2}} , \ldots, \EE^{\rho_{d}}$ be the Klyachko filtrations associated to $\EE$. Let $V_{1}, V_{2}, \ldots, V_{p}$ be the distinct one-dimensional subspaces of $E$ that appear in the Klyachko filtrations of $\EE$, i.e. each $V_{l}$ is equal to $\EE^{\rho_{j}}(i)$ for some $j \in \{1,2, \ldots, d \}$ and some $i \in {\bf{Z}}$. We now define the subsets $A_{0},A_{1},\ldots,A_{p}$ of $\{1,2, \ldots, d \}$, which intuitively classify the filtrations according to their one-dimensional subspace, as follows. For each $l \in \{ 1,2, \ldots, p \} $, we define
\[
A_{l} = \{ j \in \{1,2, \ldots, d \} \ | \ \EE^{\rho_{j}}(i) = V_{l} \text{ for some } i \in {\bf{Z}} \},
\]
and we also define
\[
A_{0} = \{1,2, \ldots, d \} \smallsetminus \bigcup_{1 \leq l \leq p} \ A_{l}.
\]
For each $j \in \{1,2, \ldots, d \}$ let $l(j)$ be the unique element of $\{0,1,2, \ldots, p \}$ such that $j \in A_{l(j)}$.
We next introduce the following collection of possibly empty subsets of $\{1,2, \ldots, d \}$:
\[
\mathcal{J} = \{ J \ | \   J \subseteq \bigcup_{1 \leq l \leq p} \ A_{l}, \text{ and $J \cap A_{l}$ has at most one element, for each $l= 1,2,\ldots,p$}   \}.
\]
Note that $\mathcal{J}$ is a finite set.
For each $j \in \{1, 2, \ldots, d \}$, we define the integers $a_{j}$ and $b_{j}$ by $a_{j} = \operatorname{max}\{ i \in {\bf{Z}} \ | \ \EE^{\rho_{j}}(i) = E \}$ and
\[
b_{j} = 
\begin{cases}
\operatorname{max}\{ i \in {\bf{Z}} \ | \ \operatorname{dim}_{k}\EE^{\rho_{j}}(i) = 1 \}  & \text{ if } j \in  \bigcup_{1 \leq l \leq p} \ A_{l} \\
a_{j} + 1 & \text{ if } j \in A_{0}.
\end{cases}
\] 
To each filtration $\EE^{\rho_{j}}$ we associate the linear functional $\lambda_{j}$ on $M_{{\bf{R}}} \times {\bf{R}}^{d+1}$ defined by
\[
\begin{matrix}
\lambda_{j} \colon & M_{{\bf{R}}} \times {\bf{R}}^{d+1} & \longrightarrow & {\bf{R}} \\
&(u,m,m_{1},\ldots,m_{d}) & \longmapsto & \frac{\langle u , v_{j} \rangle -a_{j} m - m_{j}}{b_{j}-a_{j}}
\end{matrix}
\]
for any $u \in M_{{\bf{R}}}$ and any $m,m_{1},m_{2},\ldots,m_{d} \in {\bf{R}}$.

%%%%%%%%%%%%%%%%%%%%%%%%%%%%%%%%%%%%%%%%%%%%%%%%%%%%%%%%%%%%%%%%%%%%%%%%%%%%%%%%%%%%%%%%

% Remark filtrations

\begin{remark} \label{example.concrete.klyachko}
To motivate the preceding definitions, we note that for each $j \in \{ 1 , \linebreak[0] 2 , \linebreak[0] \ldots , \linebreak[0] d \}$, and for each $m \in {\bf{Z}}_{\geq 0}$, $\om \in {\bf{Z}}^{d}$ and $u \in M$, we have
\[ {((Sym^{m}{\mathcal{E}}) \otimes \mathcal{L}^{\om})}^{\rho_{j}}(\langle u , v_{j} \rangle) = Sym^{m}_{E} \bigl( { (\mathcal{E}^{\rho_{j}}(b_{j})  )}^{ \operatorname{max} \{ 0 , \lceil \lambda_{j}(u,m,\om) \rceil \} } \bigr) \otimes L_{\om},
\]
where $\lceil x \rceil = \operatorname{min}\{ l \in {\bf{Z}} \ | \ l \geq x \}$ for any $x \in {\bf{R}}$. This equality is a reformulation of Example 3.6 in \cite{Jose1}, and it can also be verified through direct computation.
\end{remark}

%%%%%%%%%%%%%%%%%%%%%%%%%%%%%%%%%%%%%%%%%%%%%%%%%%%%%%%%%%%%%%%%%%%%%%%%%%%%%%%%%%%%%%%%

% Definition of the cone Q_{J}

We now define a rational polyhedral cone $Q_{J}$ in $M_{{\bf{R}}} \times {\bf{R}}^{d+1}$ for each $J \in \mathcal{J}$, as follows. Given such a set $J = \{j_{1},j_{2},\ldots,j_{q} \}$, the cone $Q_{J}$ is defined as the set of elements $(x,w,w_{1},\ldots,w_{d})$ satisfying the linear inequalities
\begin{align} 
&w \geq 0, \label{inequality1} \\
&\sum_{h=1}^{q}\lambda_{j_{h}}(x,w,w_{1},\ldots,w_{d}) \leq w, \label{inequality5} \\
&\lambda_{j}(x,w,w_{1},\ldots,w_{d}) \leq 0 \text{ for each } j \in \{1,2,\ldots,d \} \smallsetminus \bigcup_{1\leq h \leq q} A_{l(j_{h})}, \label{inequality2} \\
&\lambda_{j_{h}}(x,w,w_{1},\ldots,w_{d}) \geq 0 \text{ for each } h \in \{ 1,2, \ldots, q \},  \label{inequality3} \\
&\lambda_{j_{h}}(x,w,w_{1},\ldots,w_{d}) \geq \lambda_{j}(x,w,w_{1},\ldots,w_{d})  \text{ for each } h \in \{1,2, \ldots, q \} \label{inequality4}  \\
& \hspace{80 mm} \text{ and each } j \in A_{l(j_{h})}, \notag 
\end{align}
where $(x,w,w_{1},\ldots,w_{d})$ are the coordinates in $M_{{\bf{R}}} \times {\bf{R}}^{d+1}$.
The following lemma motivates the definition of these cones.

%%%%%%%%%%%%%%%%%%%%%%%%%%%%%%%%%%%%%%%%%%%%%%%%%%%%%%%%%%%%%%%%%%%%%%%%%%%%%%%%%%%%%%%%%

% Lemma H Q N

\begin{lemma} \label{lemma}
For each $g = (u,m,\om) \in M \times {\bf{Z}} \times {\bf{Z}}^{d}$ we have:
	\item[{\bf{(a)}}] If $(u,m,\om)$ belongs to the cone $Q_{J}$, for some $J = \{j_{1},j_{2},\ldots,j_{q} \} \in \mathcal{J}$, then
\[
ev_{t_{0}} \left( H^{0}(X,Sym^{m}\EE \otimes \LL^{\om})_{u} \right) = Sym^{m}_{E}\left(V_{l(j_{1})}^{\left\lceil \lambda_{j_{1}}(g) \right\rceil},V_{l(j_{2})}^{\left\lceil \lambda_{j_{2}}(g) \right\rceil},\ldots,V_{l(j_{q})}^{\left\lceil \lambda_{j_{q}}(g) \right\rceil}\right) \otimes \ L_{\om}.
\] % Lemma H
	\item[{\bf{(b)}}] If $H^{0}(X,Sym^{m}\EE \otimes \LL^{\om})_{u} \neq 0$, then $(u,m,\om)$ belongs to the cone $Q_{J}$ for some $J \in \mathcal{J}$. % Lemma Q
	\item[{\bf{(c)}}] Assume that $(u,m,\om) \in c \, (M \times {\bf{Z}} \times {\bf{Z}}^{d})$ belongs to the cone $Q_{J}$, for some $J = \{j_{1}, \linebreak[0] j_{2}, \linebreak[0] \ldots,j_{q} \} \in \mathcal{J}$, where $c = \operatorname{lcm} \{ b_{j}-a_{j} \, | \, j=1,2,\ldots,d \}$. Then $H^{0}(X,Sym^{m}\EE \otimes \LL^{\om})_{u} \neq 0$. % Lemma N
\end{lemma}
\begin{proof}
	\item[{\bf{(a)}}] From the definition of the Klyachko filtrations of $Sym^{m}\EE \otimes \LL^{\om}$, we have that
\[
ev_{t_{0}} \left( H^{0}(X,Sym^{m}\EE \otimes \LL^{\om})_{u} \right) = \bigcap_{1 \leq j \leq d} (Sym^{m}\EE \otimes \LL^{\om})^{\rho_{j}}(\langle u, v_{j}\rangle).
\]
For each $j \in \{1,2,\ldots,d \} \smallsetminus \bigcup_{1\leq h \leq q} A_{l(j_{h})}$, the inequality $\lambda_{j}(g) \leq 0$ implies that $(Sym^{m}\EE \otimes \LL^{\om})^{\rho_{j}}(\langle u, v_{j}\rangle) = Sym^{m}E \otimes L_{\om}$. For each $h \in \{1,2,\ldots,q\}$, the inequality $\lambda_{j_{h}}(g) \geq 0$ implies that $(Sym^{m}\EE \otimes \LL^{\om})^{\rho_{j_{h}}}(\langle u, v_{j_{h}}\rangle) = Sym^{m}_{E} (V_{l(j_{h})}^{\left\lceil \lambda_{j_{h}}(g) \right\rceil}) \otimes L_{\om}$. Similarly, for each $h \in \{1,2, \ldots, q \}$ and  each $j \in A_{l(j_{h})}$ the inequality $\lambda_{j_{h}}(g) \geq \lambda_{j}(g)$ implies that 
\[
(Sym^{m}\EE \otimes \LL^{\om})^{\rho_{j}}(\langle u, v_{j}\rangle) \supseteq (Sym^{m}\EE \otimes \LL^{\om})^{\rho_{j_{h}}}(\langle u, v_{j_{h}}\rangle).
\]
Therefore,
\begin{equation*}
\begin{split}
& ev_{t_{0}} \left( H^{0}(X,Sym^{m}\EE \otimes \LL^{\om})_{u} \right) \\
&\hspace{12 mm} = \bigcap_{0 \leq l \leq p} \bigcap_{j \in A_{l}} (Sym^{m}\EE \otimes \LL^{\om})^{\rho_{j}}(\langle u, v_{j}\rangle) 
=\bigcap_{1 \leq h \leq q} \bigcap_{j \in A_{l(j_{h})}} (Sym^{m}\EE \otimes \LL^{\om})^{\rho_{j}}(\langle u, v_{j}\rangle) \\
&\hspace{27 mm} = \bigcap_{1 \leq h \leq q} (Sym^{m}\EE \otimes \LL^{\om})^{\rho_{j_{h}}}(\langle u, v_{j_{h}}\rangle) 
= \bigcap_{1 \leq h \leq q} Sym^{m}_{E}(V_{l(j_{h})}^{\left\lceil \lambda_{j_{h}}(g) \right\rceil}) \otimes L_{\om} \\
&\hspace{42 mm}= Sym^{m}_{E}\left(V_{l(j_{1})}^{\left\lceil \lambda_{j_{1}}(g) \right\rceil},V_{l(j_{2})}^{\left\lceil \lambda_{j_{2}}(g) \right\rceil},\ldots,V_{l(j_{q})}^{\left\lceil \lambda_{j_{q}}(g) \right\rceil}\right) \otimes L_{\om}.
\end{split}
\end{equation*}
	\item[{\bf{(b)}}] We define
\[
I =  {\big\{} \, l \in \{1,2,\ldots,p\}  \ | \   \operatorname{max} \, \{ \lambda_{j}(u,m,\om)  \ | \  j \in A_{l} \} \geq 0 \, {\big\}}.
\]
We can assume that $I \neq \phi$, since otherwise $(u,m,\om)$ belongs to $Q_{J}$ for $J = \phi \in \mathcal{J}$. Let $l_{1},l_{2},\ldots,l_{q}$ be the distinct elements of $I$, and for each $h \in \{1,2,\ldots,q\}$ let us choose $j_{h} \in A_{l_{h}}$ such that $\lambda_{j_{h}}(u,m,\om) = \operatorname{max}\{ \lambda_{j}(u,m,\om)\  \ | \  j \in A_{l_{h}} \}$. Clearly, the set $J =_{def}\{ j_{h}  \ | \  h = 1,2,\ldots, q  \}$ belongs to $\mathcal{J}$. Then $(u,m,\om)$ satisfies (\ref{inequality1}) and (\ref{inequality5}) since \[
\bigcap_{1 \leq h \leq q} Sym^{m}_{E}\bigl(V_{l(j_{h})}^{\left\lceil \lambda_{j_{h}}(u,m,\om) \right\rceil} \bigr) \otimes L_{\om} = ev_{t_{0}} \left( H^{0}(X,Sym^{m}\EE \otimes \LL^{\om})_{u} \right) \neq 0,
\]
and it satisfies (\ref{inequality2})-(\ref{inequality4}) by the definitions of $I$ and $J$. Thus, $(u,m,\om) \in Q_{J}$.
	\item[{\bf{(c)}}] By {{(a)}}, it suffices to show that $Sym^{m}_{E}\bigl(V_{l(j_{1})}^{ \lambda_{j_{1}}(u,m,\om) },V_{l(j_{2})}^{ \lambda_{j_{2}}(u,m,\om) },\ldots,V_{l(j_{q})}^{ \lambda_{j_{q}}(u,m,\om) }\bigr)$ is non\-zero, which is true by (\ref{inequality5}) and Lemma \ref{lemma:M} (b).
\end{proof}

%%%%%%%%%%%%%%%%%%%%%%%%%%%%%%%%%%%%%%%%%%%%%%%%%%%%%%%%%%%%%%%%%%%%%%%%%%%%%%%%%%%%%%%%%

%\vspace{6 mm}

%%%%%%%%%%%%%%%%%%%%%%%%%%%%%%%%%%%%%%%%%%%%%%%%%%%%%%%%%%
%%       Finite generation of the Cox ring of \pe       %%  
%%%%%%%%%%%%%%%%%%%%%%%%%%%%%%%%%%%%%%%%%%%%%%%%%%%%%%%%%%

\section{Finite generation of the Cox ring of $\pe$}

In the next theorem we prove that any Cox ring of $\pe$, in the sense of Definition \ref{definition.cox.ring}, is finitely generated. As a corollary we obtain that $\pe$ is a Mori dream space as defined by Hu and Keel in \cite{HK}, if the toric variety $X$ is projective and $\Delta$ is simplicial (i.e. each cone in $\Delta$ is spanned by as many vectors as its dimension). Throughout this section we use the notation from \S \ref{TVBCR}-\ref{lemmas}.

%%%%%%%%%%%%%%%%%%%%%%%%%%%%%%%%%%%%%%%%%%%%%%%%%%%%%%%%%%%%%%%%%%%%%%%%%%%%%%%%%%%%%%%%%

% Theorem VERONESE PROOF.

\begin{theorem} \label{theorem}
Any Cox ring of the projectivization $\pe$ of a rank two toric vector bundle $\EE$ over an arbitrary toric variety $X$ is finitely generated as a $k$-algebra.
\end{theorem}
\begin{proof}
By Remark \ref{smooth} we can assume that $X$ is smooth. It suffices to find a finite set of generators for the $k$-algebra $R^{(c)}$ from (\ref{veronese.cox.ring}), where $c = \operatorname{lcm} \{ b_{j}-a_{j} \, | \, j=1,2,\ldots,d \}$. 
For each set $J \in \mathcal{J}$, let $G_{J} \subseteq c \, (M \times {\bf{Z}}^{d+1})$ be a finite set of generators for the semigroup $Q_{J} \, \cap \, c \, (M \times {\bf{Z}}^{d+1})$. For each $g \in M \times {\bf{Z}}^{d+1}$, let $\beta_{g}$ be a $k$-basis of $R_{g}$. We claim that the finite set
\[
\beta =_{def} \bigcup_{J \in \mathcal{J}} \ \bigcup_{g \in G_{J}} \beta_{g}
\]
generates $R^{(c)}$ as a $k$-algebra.
In order to prove the claim, consider $(u,m,\om) \in c \, (M \times {\bf{Z}} \times {\bf{Z}}^{d})$ such that $H^{0}(X,Sym^{m}\EE \otimes \LL^{\om})_{u} \neq 0$. By Lemma \ref{lemma} (a) and (b), there exist $J = \{ j_{1},j_{2}, \ldots, j_{q} \} \in \mathcal{J}$, such that $(u,m,\om) \in Q_{J}$ and
\[
ev_{t_{0}} \left( H^{0}(X,Sym^{m}\EE \otimes \LL^{\om})_{u} \right) = Sym^{m}_{E}\left( V_{l(j_{1})}^{ \lambda_{j_{1}}(u,m,\om) }  ,V_{l(j_{2})}^{ \lambda_{j_{2}}(u,m,\om) } ,\ldots,V_{l(j_{q})}^{ \lambda_{j_{q}}(u,m,\om) }\right) \otimes \ L_{\om}.
\]
Now, fix an expression $(u,m,\om) =  \sum_{g \in G_{J}} c_{g} g$, where $c_{g} \in {\bf{Z}}_{\geq 0}$ for each $g \in G_{J}$. Let $g = (u_{g},m_{g},\om_{g})$ for each $g \in G_{J}$, be the corresponding coordinates in $M \times {\bf{Z}} \times {\bf{Z}}^{d}$. From Lemma \ref{lemma} (a) applied in the cone $Q_{J}$, we get that for each $g \in G_{J}$,
\begin{align*} 
ev_{t_{0}} \bigl( H^{0}(X,Sym^{m_{g}}\EE \otimes \LL^{\om_{g}})_{u_{g}} \bigr) &= Sym^{m_{g}}_{E}\left(V_{l(j_{1})}^{ \lambda_{j_{1}}(g) },V_{l(j_{2})}^{ \lambda_{j_{2}}(g) },\ldots,V_{l(j_{q})}^{ \lambda_{j_{q}}(g) }\right) \otimes L_{\om_{g}}. 
\end{align*}
Since $\sum_{g \in G_{J}} c_{g} \lambda_{j_{h}}(g) = \lambda_{j_{h}}(u,m,\om)$ for each $h \in \{ 1,2,\ldots, q \}$, it follows by Lemma \ref{lemma} (c) and  Lemma \ref{lemma:M} that in the commutative diagram
\vspace{-8 mm}
\begin{center}
% \ar@{shape}[dir]  Shapes: |-> , {-}, {^{(}->}, etc. Directions: u,d,l,r.
\[
\entrymodifiers={+!!<0pt,\fontdimen22\textfont2>} 
\xymatrixcolsep{5pc}
\xymatrix{
% A
\bigotimes\limits_{g \in G_{J}}{(R_{g})^{\otimes c_{g}}} 
\ \ar[d]^{} 
\ar[r]^-{} & \ 
% B
R_{(u,m,\om)} 
\ar@{^{(}->}[d]^-{ev_{t_{0}}} \\
% C
\bigotimes\limits_{g \in G_{J}} (Sym^{m_{g}}E \otimes L_{\om_{g}})^{\otimes c_{g} } \ 
\ar[r]^-{} & \ 
% D
Sym^{m}E \otimes L_{\om} }
\] %\vspace{-8 mm }
\end{center}
the images of $\bigotimes_{g \in G_{J}}{(R_{g})^{\otimes c_{g}}}$ and $R_{(u,m,\om)}$ in $Sym^{m}E \otimes L_{\om}$ coincide. The injectivity of $ev_{t_{0}}$ implies that $\bigotimes_{g \in G_{J}}{(R_{g})^{\otimes  c_{g} }}$ surjects onto $R_{(u,m,\om)}$,
%\[
%H^{0}(X,Sym^{m}\EE \otimes \LL^{\om})_{u} = R_{(u,m,\om)} = \bigl( \bigotimes\limits_{g \in G_{J}}{(R_{g})^{\otimes \left\lfloor %c_{g} \right\rfloor}}\bigr) \otimes R_{g'} \subseteq k[\beta],
%\]
and this completes the proof.
\end{proof}

\begin{corollary}
The projectivization $\pe$ of a rank two toric vector bundle $\EE$ over the projective simplicial toric variety $X$ is a Mori dream space.
\end{corollary}
\begin{proof}
The finite generation of any Cox ring of $\pe$ in the sense of Hu and Keel is a consequence of Theorem \ref{theorem}. Any simplicial toric variety is ${\bf{Q}}$-factorial, and a projective bundle over a ${\bf{Q}}$-factorial variety is again ${\bf{Q}}$-factorial, hence the additional conditions follow at once from the hypotheses.
\end{proof}

The aim of Theorem \ref{theorem}, and of our work in this article, is to give further steps toward investigating the finite generation of the Cox rings of projectivizations of toric vector bundles of arbitrary rank (see Question 7.2 in \cite{PTVB}). In the articles \cite{Knop} and \cite{Suess}, the authors prove the finite generation of the Cox rings of $T$-varieties of complexity one. In particular, they provide alternative proofs of Theorem \ref{theorem}, since the projectivization of a rank two vector bundle has complexity one as a $T$-variety. We expect that our proof will provide new insights for treating this question in the case of toric vector bundles of arbitrary rank.

%\vspace{6 mm}

%%%%%%%%%%%% References %%%%%%%%%%%%%
%%
%<Author name> is written as Initial of Given Name, and Family Name.
%<Title> is written in roman letters.
%<Journal name> should be abbreviated according to
% the MR Serials Abbreviations List of Mathematical Reviews:
% (Abbreviations of Names of Serials; Mathematical Reviews).
%For <pages>, use en-dash "--" between page numbers.
%%

\bigskip

\end{document}